\documentclass[10pt]{amsart}
\usepackage[utf8]{inputenc}
\usepackage[english,russian]{babel}

\usepackage{amsmath}
\usepackage{amssymb}
\usepackage{amsfonts}

\newtheorem{theorem}{Theorem}

\begin{document}
\renewcommand{\refname}{References}
\begin{center}
\textbf {ON SOME NEW SHARP RESULTS ON TOEPLITS OPERATORS IN THE UNIT DISK}\\

\vspace{8pt}

\textbf{R.F. Shamoyan, V.A. Bednazh}\\
\end{center}

\vspace{8pt}

\noindent{\small {\bf Abstract.} We provide new sharp results on the action of Toeplitz operators from Triebel and Besov spaces to new BMOA-type function spaces on the unit disk.  In this paper we consider $s\geq1$ case in previous papers $s<1$ was covered for $BMOA_{s,q}$ and $BMOA_{s}^p$ spaces. We modify а little our previously known proofs.

\vspace{10pt}

\noindent{\small {\bf Keywords:} Toeplitz operators, Besov type spaces, Lizorkin-Triebel type spaces, analytic functions.}

\section{INTRODUCTION}

In this note we will extend our previously known sharp theorems on the action of Toeplits operators into BMOA type function spaces in the unit disk (we consider $s\geq 1$ here). More precisely we provide new sharp results on the action of Toeplits operators from mixed norm analytic  function spaces into new BMOA type classes in the unit disk.
For that reason we modify the previously known proof, provided earlier by first author
in classical  function spaces.

Proofs of our sharp results and proofs of \cite{Romi 2} are based mainly on similar type ideas.

We introduce new BMOA type spaces in the  unit disk as follows.

$$
BMOA _{s,q}(U)= \Bigg\{  f \in H^s(U):||f||_{BMOA_{s,q}}=
$$
$$=\sup\limits_{z\in U} \left(\int_{T}\frac{|f(\xi)-f(z)|^{s}}{|1-\xi\bar{z}|^{q}} dm(\xi)(1-|z|^2)\right)^{1/s}, 0<q<\infty, 1 \leq s<\infty  \Bigg\};
$$

$$
BMOA _{s}^p(U)= \Bigg\{  f \in H^s(U):||f||_{BMOA_{s}^p}=
$$
$$=\sup\limits_{z\in U} \left(\int_{T}\frac{|f(\xi)-f(z)|^{s}}{|1-\xi\bar{z}|^{2}} dm(\xi)(1-|z|^2)^p\right)^{1/s}, 0<p<\infty, 1\leq s<\infty  \Bigg\}
$$

(see definitions of all objects  below). 
It is easy to see that in particular values 
of parameters quazinorms of these analytic  spaces in the unit disk  coincide
with the so -called Garsia norm in BMOA (see \cite{2}}, \cite{3}- \cite{Romidisertacija}).

The intention of this short paper to show new sharp results on the action of  $T_{\varphi}$ Toeplitz operators in  some new BMOA type spaces in the unit disk. We provide a necessary and sufficient
condition. on the symbol of  $T_{\varphi}$ operator. Note such type results have various applications. 
Various results on BMOA type function spaces can be seen in papers   \cite{2}, \cite{3}- \cite{Romidisertacija}. Various results on the action of  $T_{\varphi}$ Toeplitz operators can been seen in recent papers \cite{Alexandrov},  \cite{2},  \cite{3},   \cite{Romidisertacija}  on various 
new and  classical analytic function spaces in the unit disk. We refer to \cite{2},  \cite{HarutyunyanShamoyan}, \cite{3}  and  \cite{Romidisertacija}  for some applications of such type results in analytic function spaces.

 Let further $ U =\left\lbrace z \in C, |z| < 1\right\rbrace  $ or $D$ be the unit disk on a complex plane C, T be the unit circle on C.  Let also further $ I=(0,1)$.
Let further H(U) be the space of all analytic functions in U .

 In this paper we as usual denote by $ \mathcal{D}^{\alpha} $ for any real $ \alpha $ the fractional derivative of analytic $ f $ function in the unit disk,
	\begin{gather*}
		\mathcal{D}^{\alpha} f(z) = \sum\limits_{k = 0}^{\infty} (k+1)^{\alpha} a_{k} z^{k},  z \in U  \\
	\text{for any analytic $ f $ function}, f(z) = \sum\limits_{k = 0}^{\infty} a_{k} z^{k},  \alpha >-1, \alpha\in \mathbb{R} \quad (see  [5]).
	\end{gather*}

 Note if $f\in H(U)$ then for any $s\in \mathbb{R}$, \ $D^s f\in H(U)$.

	We define Lusin cone in a usual manner as follows (see  \cite{3}, \cite{OrtegaFabrega}).
	\begin{gather*}
		\Gamma_{\alpha} (\xi) = \left\lbrace z \in U,\: \left| 1 - z \xi\right|  < \alpha (1 - \left| z\right| ) \right\rbrace, \text{where}  \;\alpha > 1, \xi \in T.
	\end{gather*}
	
	 We refer to   \cite{Alexandrov},  \cite{3}, \cite{OrtegaFabrega}  for further details on this object.
	
	The Hardy spaces, denoted by $H^p(U) \; (0 < p \leq \infty)$, are defined as usual (see  \cite{Romi}) by 
	$$H^p(U) = \left\lbrace f \in H(U) \colon  \sup\limits_{0 < r < 1} M_p (f, r) < \infty  \right\rbrace, $$
	
	where $$M^p_p (f, r) = \int_{T} \left| f(r \xi)\right|^p dm_1 (\xi), \;\; M_{\infty} (f, r) = \max\limits_{\xi \in T} \left| f (r \xi)\right|, r \in (0, 1), f \in H(U).$$

	For $\alpha > -1, 0 < p < \infty,$ recall that the weighted Bergman space $A^p_{\alpha} (U)$ consists of all holomorphic functions on the unit disk satisfying the condition
	\begin{gather*}
		\left\| f \right\|^p_{A^p_{\alpha}} = \int_{U} \left| f(z) \right|^p (1 - \left|  z \right|^2 )^{\alpha} dm_{2}(z) < \infty \;\; \text{(see [5, 7, 8, 9])}.
	\end{gather*}

Let further $H(U)$ be the space of all analytic functions in $U$. Let further also (see  \cite{OrtegaFabrega},  \cite{4})

$$F_\alpha^{p,q}(U)=
\left\{f \in H(U):\|f\|^{p}_{F_{\alpha}^{p,q}}=\int_{T}\left(\int_{I}|D^{m}f(r\xi)|^q(1-r)^{(m-\alpha) q-1}dr\right)^{\frac pq}d\xi<\infty\right\},$$

where $0<p, \ q<\infty$, $m> \alpha, \alpha \in \mathbb{R}$,
be the holomorphic Lizorkin-Triebel space, (see, for example, \cite{OrtegaFabrega},  \cite{4}).

Let
  $$
  F_{\alpha,k}^{p,q}(U)=\left\{f\in H(U):\|D^k f\|_{F_\alpha^{p,q}}<\infty\right\},
0<p,q,\alpha<\infty, k\in \mathbb{N}.$$

Note that we can easily show $F^{p,q}_\alpha$ general mixed norm analytic function
spaces in the unit disk are Banach spaces for all values of $p$ and $q,$ if
$\min(p,q)>1$ and they are complete metric spaces for all other values of $p$ and $q$.

 Note (see \cite{Alexandrov},  \cite{2}, \cite{OrtegaFabrega},  \cite{4})  for particular case $p=q$ we have Bergman classical class,
for $q=2$ we have so-called Hardy-Lizorkin space $H_\beta^p$ for some  $\beta$ that is,
$H_\beta^p=\{f\in H(U): D^\beta f\in H^p\}$, $0<p\leq\infty$, $\beta>0$, where $D^\beta$ is a
fractional derivative of analytic $f$ function in $U$. Note (see definitions bellow)
for this particular cases the action of $T\varphi$ classical Toeplitz operator is well-studied
in unit disk, unit ball, unit polydisk and unit disk. We study $T_{\varphi}$ operators in more general
$F_\alpha^{p,q}$ type spaces in the unit disk. Our main sharp result provide some criteria for symbol of $T\varphi$
to obtain boundednes of $T_{\varphi}$ in mentioned type analytic spaces.


Various sharp results on action of Teoplitz and other  operators can be seen in papers of various authors in
various functional spaces in the unit ball, polydisk and unit disk. We mention, for example,
the following papers \cite{HarutyunyanShamoyan} and  \cite{Romidisertacija},  where such type sharp results can be seen for various
cases of $F_\alpha^{p,q}$ spaces namely in Bergman type and in Hardy type spaces in the unit ball, polydisk 
and in the unit disk. We also note similar type results in for particular values of parameters are well-known also in other domains
(see, for example,  \cite{Romidisertacija}).

Such type sharp result on boundedness of Toeplitz operators also have various applications
(see, for example \cite{HarutyunyanShamoyan}, \cite{Romidisertacija}).

We remind the reader the standard definition of Toeplitz $T_{h}$ operators in the unit disk.

Let $h\in L^1(T)$. Then we define Toeplitz $T_{h}$ operator as an integral operator

$$(T_{h} f)(z)=\frac{1}{(2\pi)}\int_{T}\frac{f(\xi)h(\xi)}{(1-\bar{\xi} z)}dm(\xi),$$
$z\in U$.

We stress that behavior of the operators in the unit polydisk is substantially different
from the action of $T\varphi$ operators in the unit ball in $\mathbb{C}^n$ (see 
\cite{HarutyunyanShamoyan} for example).
Our intention to set criteria for the action of Toeplitz $T_{\varphi}$ operators
from $F_{\alpha,k}^{p,q}(U)$ into BMOA type
spaces in the unit disk, under the assumption that $\varphi$ is holomorphic, $\varphi\in H(U)$ (with some restriction
on symbol of Toeplitz operator).

We define some new function spaces in the unit disk for formulation of our main result in the unit disk.

$$
A_{\alpha,m}^s(U)=
\left\{f\in H(U):\|f\|_{A_{\alpha,m}}^s=\int_{U}\left|(D^m f)(z)\right|^s
(1-|z|)^{\alpha-1}dm_{2}(z)<\infty\right\},
$$
$m\in \mathbb{N}, \ 0<s, \alpha<\infty$ (Bergman-Sobolev space). Let further
$$H_{m}^s(U)=\left\{f\in H(U):\|D^m f\|_{H^s}<\infty, \ m\in \mathbb{R}, \ 0<s<\infty\right\}$$
be analytic Hardy-Lizorkin space in the unit disk $U$.

We denote by $B(z,r)$ the Bergman ball in $U$ (see \cite{ShH}, \cite{1}).

Note it can easily shown that these both scales of analytic function spaces in the unit disk are
Banach spaces for all values of $s, \ s\geq 1$ and they are complete metric spaces for other values of $s, \ s>0$.

These known spaces are particular cases of larger mixed  norm spaces $F_{\alpha, k}^{p,q}$ which we consider in this paper.

Throughout the paper, we write $C$ or $c$ (with or without lower indexes) to
denote a positive constant which might be different at each
occurrence (even in a chain of inequalities), but is independent
of the functions or variables being discussed.

We pay  special attention to places where different arguments from those we see in \cite{Romi 2} are needed.
 
 In \cite{Romi 2} we considered $s<1$ case, here we show same theorems for $s\geq 1$.
 
\section{MAIN RESULTS}  
 
 In this section we formulate our main results.

  \begin{theorem} 
 
 Let $s\geq1, q=2s-1, \tau=2(1-1/s); s\leq2.$ Then $(T_{\bar{\varphi}})$ operator is a bounded operator from $B_{2(1-1/s)}^{s};$

$$
||f||_{B_{\tau}^{s}}^{s}=\int\limits_{D}|(D^kf)(z)|^s\cdot(1-|z|)^{sk+1-2s}dm_2(z)
$$

 to $BMOA_{s,q}(D)$ if and only if $\varphi\in H^\infty(D)$

  \end{theorem}
 
 Proof of theorem 1.

 If $T_\varphi$ operator is bounded then $||T_\varphi f||_{BMOA_{s,q}}\leq C||f||_{B^s_\tau}, \tau \approx 2(1-1/s)$.

 We estimate $||T_\varphi f||$ from bellow and $||f||_{B_\tau^s}$ from above.
  First we show sufficient part.
 We have that following arguments from  \cite{Romi 2}

 $$
\lim\limits_{R \to 1} \int\limits_{T}\frac{|F(R\xi)-F(Rz)|^s(1-|w|)^2}{|1-\bar{w}\xi|^q}d\xi\leq C_1||\varphi||_{H^\infty}
\int\limits_T\frac{(1-|w|)}{|1-w\xi|^{q-s}}\left(\int\limits_D\frac{|D^kf(z)|(1-|z|)^{k-1}}{|1-z\xi|\cdot|1-\bar{w}z|}dm_2(z)\right)^s;
 $$
 
 where 
 $F(R\xi)=(T_\varphi(f))(R\xi), R\in(0,1), R>R_0, R_0\in (0,1).$
 
 Then we have  obviously
  $$
\left(\int\limits_{D}\frac{|D^kf(z)|(1-|z|)^{k-1} dm_2(z)}{|1-z\bar{\xi }|\cdot|1-\bar{w}z|^s}\right)^s\leq C_2\int\limits_{D}\frac{{|D^kf(z)|^s(1-|z|)^{(k-1)s} dm_2(z)}}{|1-\bar{z}\xi|^{s}}\times (s\geq1; s\leq2).
 $$
 
$$
\times\left(\int\limits_{D}\frac{ dm_2(z)}{|1-\bar{w}z|^{s'}}\right)^{s/s'}\leq C_3\int\limits_{D}\frac{{|(D^kf)(z)|^s(1-|z|)^{(k-1)s} dm_2(z)}}{|1-\bar{z}\xi|^{s}}\cdot(1-|w|)^{-s+2(s/s')}.
 $$
 
 Hence using composition lemma 7 Fubinis theorem
 
 $$
I=\left(\lim\limits_{r\to 1}\right)\left(\int\limits_{T}|F(r\xi)-F(rw)|^{s}\cdot \frac{(1-|w|)^s \cdot dm_2(\xi)}{|1-\bar{w}\xi|^q}\right) ^{1/s}\leq
$$

$$\leq C_4\left(\int\limits_{T} \frac{(1-|w|)}{|1-\bar{w}\xi|^{q-s}}\cdot\int\limits_{D} \frac{|(D^kf)(z)|^s(1-|z|)^{s(k-1)}}{|1-z\bar{\xi}|^{s}}\cdot (1-|w|)^{-s+2(s/s')}dm_2(z)dm(\xi)\right) ^{1/s}\leq
 $$

 $$\leq C_5 \int\limits_{D} |(D^kf)(z)|^s(1-|z|)^{s(k-2)+1}dm_2(z)\leq  C_6
  \begin{cases}
    ||f||_{A^{p,q}_{\tilde\tau}}\\
   ||f||_{F^{p,q}_{\tilde{\tau}}}
  \end{cases}
 $$
 
for some $\tilde{\tau}$ (see lemmas)
 
 Lets us show the reverse in our Theorem.  We simply repeat same arguments from  \cite{Romi 2}  to get the following estimates, based on Lemmas 4,5 and on estimates
 
 $$
||(f_r)(z)||_{B_{\tau}^{s}}^{s} \leq C_7(1-r)^{\frac{2-q-s}{s}+\gamma}\approx(1-r)^{\tilde{\tau}_0}, r<1, r>r_0.
 $$
 
$$ \text{of}  \quad   (f_r)(z)=\frac{(1-r)^\gamma}{(1-rz)}; r\in(1/2,1); z\in D; \gamma>\gamma_0 \quad \text{for} \quad\tau=2(1-1/s); q=(2s-1).
 $$

 We also have (see \cite{Romi 2}) 

$$
||T_\varphi fr||_{BMOA_{q,s}}=|\varphi(r)|||f_r||_{BMOA_{q,s}}\frac{r}{2\pi},   ||f_r||_{BMOA_{q,s}}\geq C_8(1-r)^{\tilde{\tau}_0}
$$
for some $\tau_0$.

Note also (see \cite{Romi 2}) $||f_r||_{BMOA_{q,s}}\geq C_9(1-r)^{\tilde{\tau}_0}$ .

 Note that
from our arguments(proof of sufficiency) it is easy to see that Toeplits operator  acts into Hardy space $H^{s}$.
We omit here easy details.

 Theorem 1 is completly  proved now by same argument as in [11].
 
 Theorem 1 can be reformed similarly also for more general cases of $A^{p,q}_{\tau}, F^{p,q}_{\tau}$
 and Herz type spaces similarly (see remarks below).
 
 Consider now the case of $s\geq1$ but for $BMOA_{s}^{q}$ spaces. Note

  $$
(BMOA_{s}^{q})(D)=\Bigg\{f\in (H^s)(D)
:||f||_{BMOA_{s}^{q}}=\left(\sup\limits_{z
\in D}\right)\left
(\int\limits_{T} \frac{|f(\xi)-f(z)|^{s} dm(\xi)
}{|1-\bar{\xi}z|^2}\cdot(1-|z|)^q\right) ^{1/q}\Bigg\};
$$
 
 $0<q<\infty; 1\leq s<\infty.$

 We have that following the proof above, arguments from \cite{Romi 2}, and Lemmas.
 
 $$I(F)= \left(\lim\limits_{r\to 1
}\right)\left
(\int \limits_{T} |F(r \xi  )
 - F(rw)|^{s} \frac{(1-|w|^2)^{q} 
}{|1-\bar{w} \xi|^2}dm(\xi)\right) ^{1/s}\leq
$$
 
 $$\leq C_{10} \left(\int\limits_{T} \int\limits_{D} |(D^kf)(z)|^s \cdot \frac{(1-|z|)^{s(k-1)}}{|1-z\bar{\xi}|^{s}}\cdot (1-|w|)^{-s+2(s/s')}\cdot\frac{(1-|w|)^q}{|1-\bar{w}\xi|^{2-s}} dm(\xi)dm_2(z)\right) ^{1/s}\leq
 $$

 $$\leq C_{11}  \int\limits_{D} |D^kf(z)|^s(1-|z|)^{s(k-1)+1-s}dm_2(z).
 $$
 
 $$
1\leq s \leq2; \quad q=4-2s; \quad q \in [0,2].  
 $$
 
 We show the reverse now. We have (see also arguments above)
 
 $$
 ||f_r(z)||_{(BMOA)^{q}_{s}(D)}\geq C_{12}(1-r)^{\gamma+\frac{q-1-s}{s}}; \gamma>\gamma_0,
 $$
 
 and to show this we follow arguments of  \cite{Romi 2}
 
 $$
\left( ||(f_r)||_{BMOA^{q}_{s}} \right)=
 (1-r)^{\gamma}\left((\sup\limits_{ z\in D}) 
 \int\limits_{T}\left|\frac{1}{(1-rz)}-
 \frac{1}{(1-r\xi)}\right|^{s}\times\frac{(1-|z|)^q}{(1-\xi z)^2}d\xi\right)^{1/s}= 
 $$

 $$
= ((1-r)^{\gamma})(\sup\limits_{ z\in D}) 
 \int\limits_{T} \frac{(dm(\xi))\cdot(1-|z|)^q}{|1-r\xi|^s|1-\xi z|^{2-s}|1-rz|^s}\geq C_{13} (1-r)^{\frac{\gamma+q-s-1}{s}}=(1-r)^{\tau_0}, s\geq 1. 
 $$ 
 
 Using that $(\sup\limits_{z\in D})\geq (\sup\limits_{z=r})$ and lemmas   also

 $$
  ||f_r||_{B^s_{2(1-\frac{1}{s})}}\leq C_{14}(1-r)^{\tau_0}.
 $$
 
 We finish    the proof using same arguments as in [11] or in previous theorem 1.

 We give formulation

   \begin{theorem} 
 
 Let $(s)\geq1, q=4-2s, q\in[0,2]; s\leq2.$.
 Then $(T_{\bar{\varphi}})$ is a bounded operator from $B_{2(1-1/s)}^{s}$ to $BMOA_{s}^{q}(D)$ if and only if $\varphi\in H^\infty(D)$

The  formulation for mixed norm $A^{pq}_{\alpha}$ spaces is the following.

   \begin{theorem}
 
 Let $(s)\geq1, q=2s-1, \tau=2(1-1/s); s\leq2.$  
 Then $(T_{\varphi})$ is a bounded operator from $(D^kA^{p,\tilde{q}}_{\alpha})$ or $D^kF^{p\tilde{q}}_{\alpha}$ to $BMOA_{s,q}$ if and only if $\varphi\in H^\infty(D)$ , for $(max)(p,\tilde{q})\leq s<\infty, p=s/q,$ $ \alpha=k-1/s, k\geq 1/s$.
  \end{theorem}

The proof of these general result follow from $p=\tilde{q}$ case and embedding in lemmas we leave this task to readers. We add some words for Herz spaces.

Note that the following embeddings holds

$$
\int\limits_{D}|f(z)|^p(1-|z|)^\alpha dm_2(z)\leq C_1\sum\limits_{k\geq 0}\left(\int\limits_{B(a_k,r)}|f(z)|^p(1-|z|)^\alpha dm_2(z)\right)^{q/p}
$$

$
p\in(0,\infty), q\leq p, \alpha>-1;
$

$$
\int\limits_{D}|f(z)|^p(1-|z|)^\alpha dm_2(z)\leq C_2\int\limits_U \left(\int\limits_{B(v,r)}|f(z)|^q(1-|z|)^{\alpha+2} dm_2(z)\right)^{q/p}dv;
$$

$
 \alpha>-1, 0<p,q<\infty.
$
these are quazinorms of  Herz spaces.

Since similar sharp result is valid for $s\geq 1, BMOA^s_q, A^{p,q}_{\alpha}, F^{p,q}_{\alpha}$ tripple. 
This gives us a chance to easily extend our theorems 1, 2 also to Herz type function analytic spaces in the Unit disk $U$.

 Consider other spaces with following quazinorms

$$
\left(\sup\limits_{|z|<1}\right)\left|(D^k f)(z)\right|\cdot(1-|z|)^{k-s}<\infty \qquad\qquad\qquad\qquad (A)
$$

 or
 
 $$
\int\limits_{0}^{1}\left(M_{\infty}(D^kf,r)\right)^p\cdot(1-r)^{\tilde{k}-\alpha} dr <\infty; \qquad\qquad\qquad\qquad (B)
 $$

$0<p<\infty, \alpha>-1, k>s, \tilde{k}>\alpha$, these are another analytic Besov spaces,

or

$$
\left(\sup\limits_{r<1}\right)\left(\int\limits_{T}|(D^k f)(r\xi)|^pd\xi\right)^{1/p}\cdot(1-r)^{k-s}<\infty; \qquad\qquad\qquad\qquad (C)
$$

$k>s,0<p<\infty.$
  \end{theorem}

Denote them by $X$, We wish to find sharp conditions of $\varphi$, so that $(T_{\varphi})$  acting from $X$ to $(BMOA_{s,q})(D)$ or $(BMOA^q_s)(D)$ as bounded operator. We can  modify the proof we suggested above for these type of analytic Besov type spaces also.

 We turn to other function spaces (A), (B), (C). We give  partial answers in Theorem 4.

 For $s\leq 1$ case we arrived  similarly (see above) and \cite{Romi 2}.. 

 $$||T_\varphi f||_{BMOA_{s,q}}\leq C_1\int\limits_{T} \left(\int\limits_{D} \frac{ |D^kf(z)|^s (1-|z|)^{ks+s-2} (1-|w|)dm_2(z)}{|1-z\bar{\xi}|^{s}|1-\bar{w}\xi|^{s}|1-\bar{w}\xi|^{q-s}}\right) ^{1/s};
 $$

$w\in D; BMOA_{s,q}$ space.

 For $s\geq 1$ case we arrived to the estimate  (see above). 

 $$||T_\varphi f||_{BMOA_{s,q}}\leq C_2 \int\limits_{T} \left(\int\limits_{D} \frac{ |D^kf(z)|^s (1-|z|)^{s(k-1)} (1-|w|)^{2s/s'-s+1}(dm(\xi))(dm_2(z))}{|1-\xi \bar{z}|^{s}|1-\bar{w}\xi|^{s}|1-\bar{w}\xi|^{q-s}}\right) ^{1/s'}
 $$

for $ BMOA_{s,q}$ space.

Estimating further accurately we have the following chain of estimates for $s\geq 1$ case separately for $BMOA_{s,q}$ and $BMOA^q_s$ spaces. 

We have

 $$ \int\limits_{T} \int\limits_{D} \frac{ |D^kf(z)|^s (1-|z|)^{s(k-1)} (1-|w|)^{2s/s'-s+1}dm_2(z)dm(\xi)}{|1- \bar{\xi}z|^{s}|1-\bar{w}\xi|^{q-s}}\leq
 $$

$$
\leq C_3\int\limits_0^1
\left(M_{\infty}(D^kf(z))^s\cdot (1-|z|)\right)^{s(k-1)+1-s}dr= M(f);$$

$(q-s)>1; (q-s)-1=2s/s'-s+1, q=2s, s\geq1.$

The reverse implication. 

Note that for our  test function $f_r$

$$
M(f)^{\frac{1}{s}}\leq C_4(1-r)^\gamma\cdot (1-r)^{-(k+1)s+2-s+s(k-1)}=C_5(1-r)^{\gamma+\frac{2-3s}{s}}; 
$$

$$
||f_r||_{BMOA_{s,q}}\geq(1-r)^{\gamma+\frac{2-3s}{s}}.
$$

And following the proof of previous case we arrive at new sharp result for $T_\varphi$.

A version of this theorem for $(BMOA^s_q)$ spaces can be also formulated. 

 \begin{theorem}
 
 Let $s>1, q=2s, s\leq 2$.  Then $(T_\varphi)$ is a bounded operator from a Besov space with quazinorm
 
 $$
K(f)=\int\limits_{0}^{1}\left(M_{\infty}(D^kf,|z|)^s\right)\cdot(1-|z|)^{s(k-1)+(1-s)} dr <\infty,
 $$

to $(BMOA_{s,q})(D)$ if and only if $\varphi \in H^{\infty}(D)$.
 
   \end{theorem}

\section{Lemmas}

We collect some interesting facts and lemmas in this section. They are taken from [11] and important for this paper.
 
  We define analytic Besov type  and Lizorkin-Triebel type spaces in the unit disk as follows.
 
 $$A_{s}^{p,\tilde{q}}=\Big\{f\in H(U): \int\limits_{I}\left(\int\limits_T|D^kf(r
\xi)|^p d\xi\right)^{\tilde{q}/p}(1-r)^{\tilde{q}(k-s)-1}dr<\infty\Big\};$$ 
 
 $k>s, 0<p,\tilde{q}<\infty, s \in\mathbb{R}; $ and  $F_{s}^{p, \tilde{q}}$  defined similarly based on definition of $F_{\alpha, k}^{p,q}$ spaces above
 changing the order of integration  (see above).

 Indeed similar  results are valid for Herz spaces (not only $A_{\alpha}^{p,q}, F_{\alpha}^{p,q},$ spaces) For  $0<p,q<\infty, \alpha>-1$ we define analityc Herz spaces as follows. Let $B(z,r)$ be Bergman ball in $U$.  These are spaces with quazinorms. 
 
 $$\int\limits_{U}\left(\int\limits_{B(z,r)}
 |D^k f (\tilde{z})|^{p}\cdot (1-|\tilde{z}|)^{\alpha}\cdot dm_2(\tilde{z})\right)^{q/p} dm_2(z) 
 $$
 
 $$
\sum\limits_{k\geq 0}\left(\int\limits_{B(a_k,r)}
 |D^k f (\tilde{z})|^{p}\cdot (1-|\tilde{z}|)^{\alpha}dm_2(\tilde{z}) \right)^{q/p},
 $$

 where $\{a_k\}$ is an $r-$lattice in $U$.

Let
$$\tilde{A}_{\alpha}^{p,q}=\Big\{f\in H(U): \int\limits_{I}\left(\int\limits_T|f(r
\xi)|^p d\xi\right)^{q/p}(1-r)^{\alpha q-1}dr<\infty\Big\};$$ 
 
 $$\tilde{F}_{\alpha}^{p,q}=\Big\{f\in H(U): \int\limits_{T}\left(\int\limits_I|f(r
\xi)|^q (1-r)^{\alpha q-1} dr\right)^{p/q}d\xi<\infty\Big\},$$ 
$0<p,q<\infty, \alpha\in (0; \infty).$

We formulate now several lemmas which are needed for proofs of our main results.

\textbf{Lemma 1.}  (see  \cite{Romidisertacija}, \cite{Romi} for mixed norm spaces). 

Let $0<max(p,q)\leq \tilde{s}<\infty$. Then we have that 
$$\left(\int\limits_U
 |f(z)|^{\tilde{s}}\cdot (1-|z|)^{\tilde{s}(\alpha+\frac{1}{p})-2}\cdot dm_2(z)\right)^{(1/\tilde{s})}\leq C_1||f||_{\tilde{F}_{\alpha}^{p,q}}; $$
 
and
$$\int\limits_U
 |f(z)|^{\tilde{s}}\cdot (1-|z|)^{\tilde{s}(\alpha+\frac{1}{p})-2}\cdot dm_2(z)\leq C_2||f||^{\tilde{s}}_{\tilde{A}_{\alpha}^{p,q}}; $$

As corollary of Lemma 1 we have Lemma 2,3
for some positive constants $C$ and $ C_1$.

\textbf{Lemma 2.}  (see  \cite{OrtegaFabrega},  \cite{Romidisertacija}, \cite{Romi}).

Let   $0<max(p,\tilde{q})\leq \tilde{s}<\infty$;
Then we have that 

$$J=\left(\int\limits_U
 |(D^k f )(z)|^{\tilde{s}}\cdot (1-|z|)^{k \tilde{s}-q}\cdot dm_2(z)\right)^{1/\tilde{s}}\leq  C_3\left( \int\limits_{I}\left(\int\limits_T|(D^k f)(r
\xi)|^p d\xi\right)^{\tilde{q}/p}(1-r)^{\alpha \tilde{q}-1}dr\right)^{1/\tilde{q}};$$
 
 $$ \alpha=(\frac{1}{s})(k\tilde{s}-q); p=\frac{\tilde{s}}{2}; k>k_0, k_0=\frac{q}{\tilde{s}};$$
 
and 

$$J^p\leq C_4\int\limits_T\left(\int\limits_I
 |(D^k f )(r\xi)|^{\tilde{q}}\cdot (1-r)^{\alpha \tilde{q}-1}\ dr\right)^{p/\tilde{q}}d\xi.$$

 \textbf{Lemma 3.}  (see \cite{OrtegaFabrega},  \cite{Romidisertacija}, \cite{Romi}).

 Let   $0<max(p,\tilde{q})\leq \tilde{s}<\infty$;
Then we have that 

$$J=\left(\int\limits_U
 |D^k f (z)|^{\tilde{s}}\cdot (1-|z|)^{k \tilde{s}+q-3}\cdot dm_2(z)\right)^{q/\tilde{s}}\leq  C_5 | |D^k f ||_{\tilde{A}_{\alpha}^{p,\tilde{q}}}=$$
 
$$ =C_6\left( \int\limits_{0}^{1}\left(\int\limits_T|(D^k f)(r
\xi)|^p d\xi\right)^{\tilde{q}/p}(1-r)^{\alpha \tilde{q}-1}dr\right)^{1/\tilde{q}};$$
 
 $$ p=\frac{\tilde{s}}{q} ;  \alpha=k-\frac{1}{2}; k>\frac{1}{\tilde{s}};$$
 and similarly  for spaces $F_{\alpha}^{p,\tilde{q}}$.

  \textbf{Lemma 4.}  (see   \cite{OrtegaFabrega},  \cite{Romidisertacija}, \cite{Romi}).

  Let $(f_r)(z)=\frac{(1-r)^\gamma}{1-rz}, \gamma>\gamma_0, \gamma_0=\gamma_0(p,\tilde{q}, s), r\in (1/2,1)$ , then we have that 
  
  $$
||f_r(z)||_{A_{q/\tilde{s}}^{p,\tilde{q}}}\leq C_8 (1-r)^{\gamma+(2-q-s)/s};
  $$

  $$
||f_r(z)||_{F_{q/\tilde{s}}^{p,\tilde{q}}}\leq C_9(1-r)^{\gamma+(2-q-s)/s}.
  $$
 
 The proof is standard. It uses classical estimates of function theory in the unit disk.

   \textbf{Lemma 5.}  (see \cite{Romidisertacija}, \cite{Romi}).
    
    Let $F\in H(U), p\leq 1 , \beta\in(-1, \infty); w,w_1 \in U; 0\leq q, q_1<\infty,$ then we have that 
 
 $$
 \left(\int\limits_{U}\frac{|F(z)|(1-|z|)^{\beta} dm_2(z)}{|1-zw_1|^{q_1}|1-zw|^q}\right)^p \leq
 C_{10}\int\limits_{U}\frac{|F(z)|^p(1-|z|)^{\beta p+2p-2} dm_2(z)}{|1-zw_1|^{q_1 p}|1-zw|^{q p}}, w, w_1 \in U.
 $$
 
 Note in particular cases if $q=0$ or $q=1$ this result is well-known. Note this result is also well-known in the ball, polydisk. And similar  result is valid for $p>1$ (see above).
 
   \textbf{Lemma 6.}  (see  \cite{HarutyunyanShamoyan},  \cite{OrtegaFabrega}, \cite{4}, \cite{Romidisertacija}).

 Let $$(D^mf)(z)=\sum\limits_{{k\geq 0}}\frac{\Gamma(m+k+1))a_k z^k}{\Gamma(m+1)\Gamma(k+1)}$$
 
 $ m\geq 0, f \in H(U), f(z)=\sum\limits_{k\geq0}a_kz^k, z \in U.$
 
Then $D^mf\in H(U)$, if $f\in H(U), f(z)=\sum\limits_{k\geq0}a_k z^k$ and in addition  we have
 
 $$
\left(\frac{1}{2\pi}\right) \int\limits_{T}f(rt)\overline{g(rt)}dm(t)=\left(\frac{2^{m-1}}{r^2(m+1)\pi}\right)\int\limits_{0}^{r}\int\limits_{T}\left(f(r\xi)\right)\overline{D^m g(r\xi)}(1-r)^{m-1}rdrdm(\xi), m\geq0, r \in(0,1).
 $$

  \textbf{Lemma 7.} (see \cite{4}).
 
 We have the following estimates
 
 $$
\int\limits_{U}\frac{(1-|\xi|^2)^{s} dm_2(\xi)}{|1-\bar{\xi z}|^r|1-\bar{\xi}w|^v}\leq \frac{C_{11}}{|1-\bar{z}w|^{r+v-s-2}};
 $$
 
$z,w\in U, s>-1; r,v\geq0, r+v-s>2, r-s<2, v-s<2$

 $$
\int\limits_{U}\frac{(1-|\xi|^2)^s dm_2(\xi)}{|1-\bar{\xi z}|^r|1-\bar{\xi}w|^v}
\leq \frac{C_{12}}{(1-|z|^2)^
{r-s-2}|1-\xi w|^v};
 $$
 
$z,w\in U, s>-1; r,v\geq0, r+v-s>2, v-s<2< r-s;$

 $$
\int\limits_{T}\frac{d\xi}{|1-\bar{\xi z}|^s|1-\bar{\xi}w|^{q-s}}\leq \frac{C_{13}}{(|1-\bar{z}w|^{q-1}};
 $$
 
$z,w\in U, 1<q<1+s, 0<s<1; r,v\geq0, r+v-s>2, v-s<2< r-s;$

  $$
\int\limits_{T}\frac{d\xi}{|1-\xi z|^s|1-\xi w|^{q-s}}\leq \frac{C_{14}(1-|w|)^{s+1-q}}{(|1-zw|^{s}};
 $$
 
$z,w\in U, 1+2s<q<\infty, 0<s<1.$

 \end{document}